\documentclass[graybox]{svmult}

\usepackage{type1cm}        
\usepackage{makeidx}         
\usepackage{graphicx}       
                            
\usepackage{multicol}        
\usepackage[bottom]{footmisc}

\usepackage{subfig}

\usepackage{newtxtext}        
\usepackage{newtxmath}

\makeindex             
                      
\usepackage{algorithm}
\usepackage{algorithmic}
\def\PP{{{\rm l}\kern - .15em {\rm P} }}
\def\PN2{{\PP_{N}-\PP_{N-2}}}

\newcommand{\bphi}{\boldsymbol{\varphi}}

\newcommand{\btau}{\boldsymbol{\tau}}

\newcommand{\ba}{\boldsymbol{a}}

\newcommand{\bff}{\boldsymbol{f}}
\newcommand{\bFF}{{\boldsymbol F}}
\newcommand{\bG}{{\boldsymbol G}}

\newcommand{\bu}{\boldsymbol{u}}

\newcommand{\bur}{{\boldsymbol{u}}_r}

\newcommand{\bx}{\boldsymbol{x}}
\newcommand{\bX}{\boldsymbol{X}}

\newcommand{\bXr}{{\bf X}^r}

\definecolor{vargreen}{rgb}{0.0, 0.5, 0.0}

\newcommand{\deleted}[1]{{}}

\begin{document}

\title*{Reduced Order Model Closures: A Brief Tutorial} 

\author{W. Snyder, C. Mou, H. Liu, O. San, R. De Vita, and T. Iliescu}

\institute{
William Snyder 
\at Virginia Tech, 
Blacksburg, VA 24061, \email{swilli9@vt.edu}
\and 
Changhong Mou
\at University of Wisconsin-Madison, Madison, WI
53706, \email{cmou3@wisc.edu}
\and
Honghu Liu
\at Virginia Tech, 
Blacksburg, VA 24061, \email{hhliu@vt.edu}
\and 
Omer San
\at Oklahoma State University, Stillwater, OK 74078, \email{osan@okstate.edu}
\and
Raffaella De Vita
\at Virginia Tech, 
Blacksburg, VA 24061, \email{devita@vt.edu}
\and
Traian Iliescu 
\at Virginia Tech, 
Blacksburg, VA 24061, \email{iliescu@vt.edu}
}

\maketitle

\abstract{
In this paper, we present a brief tutorial on reduced order model (ROM) closures. 
First, we carefully motivate the need for ROM closure modeling in under-resolved simulations.
Then, we construct step by step the ROM closure model by extending the classical Galerkin framework to the spaces of resolved and unresolved scales.
Finally, we develop the data-driven variational multiscale ROM closure and then we test it in fluid flow simulations. 
\newline \indent
Our tutorial on ROM closures is structured as a sequence of questions and answers, and is aimed at  first year graduate students and advanced undergraduate students.
Our goal is not to explain the ``how,'' but the ``why.''
That is, we carefully explain the principles used to develop ROM closures, without focusing on particular approaches.
Furthermore, we try to keep the technical details to a minimum and describe the general ideas in broad terms while citing appropriate references for details.
}

\section{Introduction}
	\label{sec:introduction}

Reduced order models (ROMs) are computational models whose dimensions are orders of magnitude lower than the dimensions of the full order models (FOMs) (i.e., models obtained from classical numerical methods, e.g., the finite element method).
Because ROMs are relatively low-dimensional, their computational cost is orders of magnitude lower than the computational cost of FOMs.
Thus, ROMs represent a promising alternative to FOMs in computationally intensive applications, e.g., digital twins of wind farms and real time surgical procedures. 
ROMs are expected to play a key role in establishing mathematical modeling foundations for digital twins of many engineering, healthcare, and environmental systems. 
Indeed, if ROM results are nearly indistinguishable from the corresponding FOM results, then they can  
contribute as predictive tools in emerging digital twin infrastructures.
However, despite being successfully used in simple, academic test problems, ROMs have not made a significant impact in complex, practical applications.

One of the main hurdles in the ROMs' development is their notorious inaccuracy when they are used in the {\it under-resolved} regime, i.e., when the ROM's dimension (i.e., its number of degrees of freedom  (DOF)) is not large enough to capture the complex dynamics of the underlying system.
To illustrate the under-resolved regime, think of the numerical simulation of the flow around a wind farm.
This simulation with a FOM (e.g., the finite element method) generally requires millions (if not billions) of DOF.
Thus, performing shape optimization or real time control of the wind farm flow, which would require many individual FOM runs, is not feasible.
Replacing the costly FOM with a ROM would be a natural choice. 
However, in order to represent the turbulent flow dynamics in the wind farm simulation would require thousands or tens of thousands of DOF in the ROM.
Despite the ROM's cost being much lower than the FOM cost, it is still too high to allow the use of the ROM in real time control applications, where thousands of ROM runs would be required.
Thus, a practical choice would be to use much cheaper ROMs, i.e., ROMs with much fewer (e.g., hundreds or even tens) DOF.
However, these low-dimensional ROMs, although computationally efficient (and, therefore, practical), generally yield inaccurate results. 
The reason is simple: these ROMs do not have enough DOF to represent the complex dynamics of a complex flow such as the turbulent wind farm flow.

The above discussion yields the following two important conclusions:

\begin{enumerate} 
	\item The under-resolved ROM regime is critical in realistic, complex applications.
	\item  Under-resolved ROMs produce inaccurate results.
\end{enumerate}

These conclusions naturally lead to the following question:

\begin{question}{Q0}
	How do we fix the under-resolved ROMs?
\end{question}

The answer to Q0 is simple:
\begin{important}{A0}
	We develop good ROM closure models, i.e., correction terms that increase the standard ROM's accuracy.
\end{important}

To our knowledge, the first (and only) survey of ROM closure models was performed in~\cite{ahmed2021closures}, where the authors discuss dozens of ROM closures for fluids that have been developed over the last four decades.
We are not aware, however, of a tutorial on ROM closures.
This paper takes a first step at filling that gap.

This brief tutorial on ROM closures (also known as parameterizations~\cite{Berneral2017,chekroun2020variational,chekroun2015stochastic,chorin2015discrete,crommelin2008subgrid,MTV01,zanna2017scale} and hidden dynamics~\cite{pawar2020data,pawar2020evolve}) 
is structured as a sequence of simple questions and answers that lead the reader from a simple PDE to projection ROMs, and then to ROM closures. 
Our paper is aimed at first year graduate students and advanced undergraduate students.
Thus, we strive to keep the technical details to a level that is easily understood by students with a standard background in differential equations and numerical methods.
We also emphasize that our goal in this tutorial is not to explain the ``how,'' but the ``why.''
That is, we carefully explain the principles used to develop ROM closures, without focusing on particular approaches (which are carefully discussed in~\cite{ahmed2021closures}).

The rest of the paper is organized as follows:
In Section~\ref{sec:toy}, we illustrate the ROM closure modeling concept for a three-dimensional toy problem.
In Section~\ref{sec:g-rom}, we present the general algorithm used to develop the classical Galerkin ROM.
In Section~\ref{sec:closure-problem-closure-model}, we first present the ROM closure problem, and then we discuss its solution, i.e., the ROM closure model.
In Section~\ref{sec:d2-vms-rom}, we construct the data-driven variational multiscale ROM, in which available data is used to build the ROM closure model. 
In Section~\ref{sec:numerical-results}, we 
illustrate how closure modeling can significantly increase the ROM accuracy in the numerical simulation of fluid flows.
In Section~\ref{sec:mathematics}, we survey  current mathematical results for ROM closure modeling.
Finally, in Section~\ref{sec:conclusions}, we present conclusions and future research avenues.

\section{A Crash Course in ROM Closure: A Toy Problem}
    \label{sec:toy}

\begin{figure}[ht]
	\sidecaption
	\includegraphics[scale=.75]{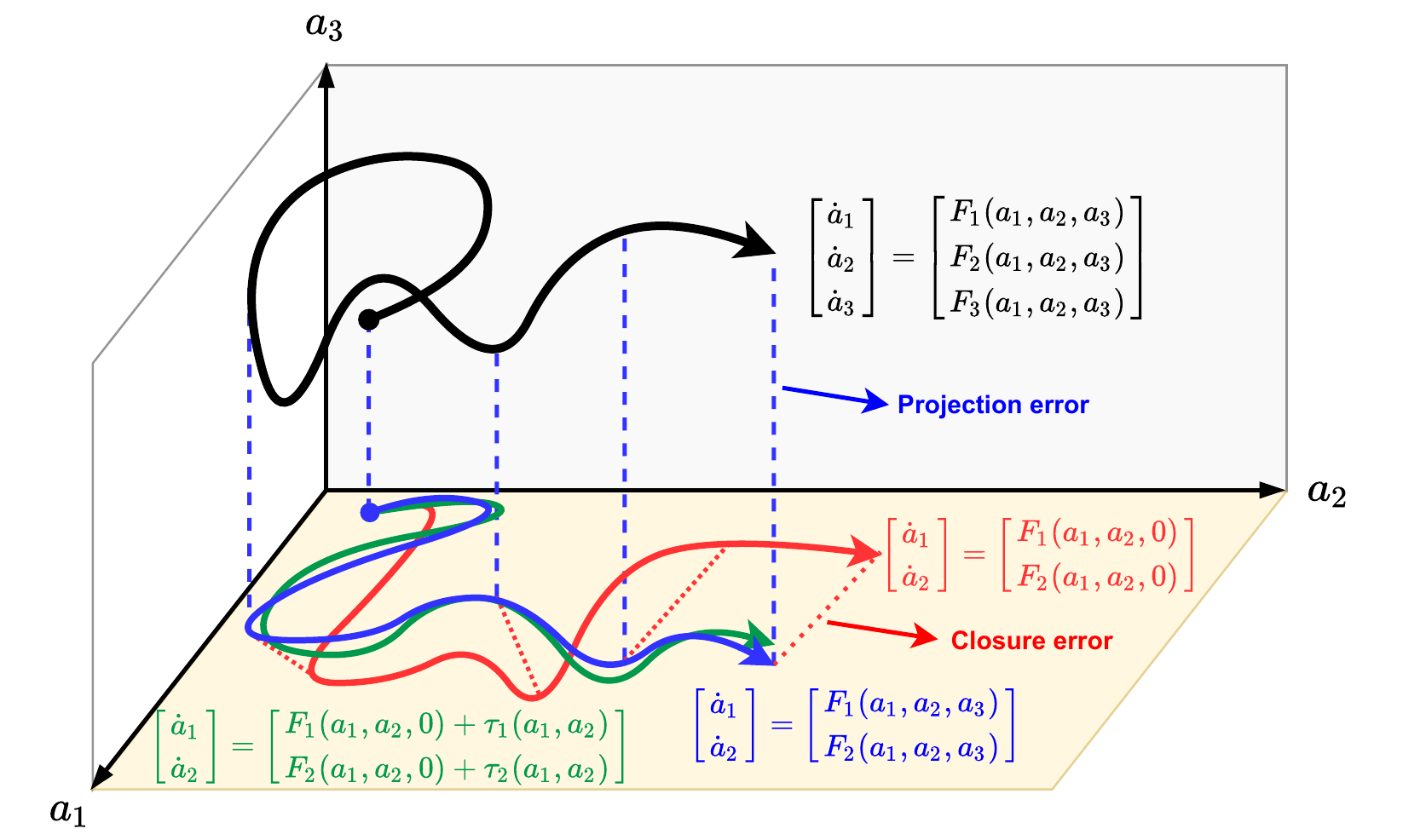}
	\caption{
	    A schematic representation of 
	    the ROM closure modeling for a three-dimensional toy problem.
	    The goal is to reduce the three-dimensional G-ROM~\eqref{eqn:schematic-2} (black curve and equations) to the most accurate two-dimensional ROM.
	    The I-ROM~\eqref{eqn:schematic-3} (blue curve and equations) is the most accurate ROM obtained in the Galerkin framework, but it is not closed (since it depends on $a_3$).
	    The two-dimensional G-ROM~\eqref{eqn:schematic-4} (red curve and equations) is closed, but it is not accurate (since we simply ignore the $a_3$ contribution).
	    The two-dimensional G-ROM supplemented with a closure model~\eqref{eqn:schematic-5} (green curve and equations) is closed and more accurate than the two-dimensional G-ROM~\eqref{eqn:schematic-4} since the closure terms $\tau_1(a_1, a_2)$ and $\tau_2(a_1, a_2)$ aim at steering the green curve toward the blue curve.
	}
	\label{fig:schematic}
\end{figure}

Before carefully presenting the ROM closure modeling in the next sections, we illustrate the underlying {\it concepts} and {\it principles} for a {\it toy problem}. 
These concepts and principles are broadly illustrated in the schematic in Fig.~\ref{fig:schematic}, which is adapted from Fig.~1 in~\cite{ahmed2022physics}.

To present our toy problem, we first assume that the FOM solution,  
$\bu^{FOM}$, can be accurately approximated by only three ROM basis functions:
\begin{eqnarray}
    \bu^{FOM}(x,t)
    \approx a_1(t) \bphi_1(\bx)
    + a_2(t) \bphi_2(\bx)
    + a_3(t) \bphi_3(\bx), 
    \label{eqn:schematic-1}
\end{eqnarray}
where $\bphi_1, \bphi_2, \bphi_3$ are the ROM basis functions, and $a_1, a_2, a_3$ are the sought time dependent coefficients.
Of course, for complex systems, one should use many more (e.g., hundreds and even thousands of) ROM basis functions to accurately approximate $\bu^{FOM}$.
However, to graphically illustrate the need for closure modeling in our toy problem, we assume that three ROM basis functions are enough.

Next, we use the three ROM basis functions in the Galerkin framework to construct the {\it Galerkin ROM (G-ROM)}.
Details regarding the G-ROM construction are given in Section~\ref{sec:g-rom}.
For the purpose of the toy problem illustration in this section, we just note that the resulting G-ROM is a three-dimensional dynamical system that can be written as follows:
\begin{equation}
\begin{aligned}
    \begin{bmatrix} \dot{a}_1\\ \dot{a}_2\\ \dot{a}_3 \end{bmatrix} = \begin{bmatrix}     F_1(a_1,a_2,a_3)\\F_2(a_1,a_2,a_3)\\F_3(a_1,a_2,a_3) \end{bmatrix},
\end{aligned} \label{eqn:schematic-2} 
\end{equation}
where $F_1, F_2,$ and $F_3$ are the components of the ROM operators, e.g., vectors, matrices, and tensors, which are presented in Section~\ref{sec:g-rom}. 
Since the three ROM basis functions yield an accurate approximation of the FOM solution in~\eqref{eqn:schematic-1}, 
the three-dimensional G-ROM in~\eqref{eqn:schematic-2} is expected to yield an accurate approximation to $\bu^{FOM}$. 
That is, solving the three-dimensional G-ROM~\eqref{eqn:schematic-2} for $a_1, a_2, a_3$, and then plugging these values back into~\eqref{eqn:schematic-1} yields an accurate approximation to $\bu^{FOM}$.  
In Fig.~\ref{fig:schematic}, the time evolution of the solution of the accurate three-dimensional G-ROM~\eqref{eqn:schematic-2} is represented as the black curve.

At this point, we invoke the need to {\it reduce the computational cost} of the three-dimensional G-ROM~\eqref{eqn:schematic-2}.
Specifically, we aim at constructing a {\it two-dimensional} ROM that is as accurate as possible (preferably, as accurate as the three-dimensional G-ROM~\eqref{eqn:schematic-2}).
For our toy problem~\eqref{eqn:schematic-1}, this amounts to constructing a dynamical system for $a_1$ and $a_2$ (assuming that the first two ROM basis functions dominate the third, as is often the case; see Section~\ref{sec:g-rom}). 

Of course, reducing the ROM dimension from three to two does not yield such a great reduction of computational time.
We emphasize, however, that we consider this reduction only to illustrate the ROM closure modeling concept for our toy problem.
In practical settings, ROMs reduce the FOM dimension by orders of magnitude.

The most natural way to construct an accurate two-dimensional ROM is to keep only the first two equations in~\eqref{eqn:schematic-2} and discard the third equation:
\begin{equation}
\begin{aligned}
    \begin{bmatrix} 
        \dot{a}_1\\ 
        \dot{a}_2 
    \end{bmatrix} 
    = \begin{bmatrix}              
        F_1(a_1,a_2,a_3)\\ 
        F_2(a_1,a_2,a_3) 
    \end{bmatrix}.
\end{aligned} 
\label{eqn:schematic-3} 
\end{equation}
Mathematically, this amounts to first using a Galerkin expansion for all three ROM basis functions (i.e., using~\eqref{eqn:schematic-1}), and then using a Galerkin projection onto only the first two basis functions (instead of projecting onto all three basis functions, as done in~\eqref{eqn:schematic-2}).

In Fig.~\ref{fig:schematic}, the time evolution of the solution of the efficient, two-dimensional ROM~\eqref{eqn:schematic-3} is represented as the blue curve.
Of course, since we perform a Galerkin projection only onto the first two basis functions, we incur an error, which we denote as the (Galerkin) {\it projection error} (the blue dashed lines in Fig.~\ref{fig:schematic}).
Nevertheless, it stands to reason that, in the Galerkin framework with the basis $\{ \bphi_1, \bphi_2, \bphi_3 \}$, the two-dimensional ROM~\eqref{eqn:schematic-3} is the most accurate two-dimensional ROM we can hope to get. 
This is why we call the two-dimensional ROM~\eqref{eqn:schematic-3} the {\it ideal ROM (I-ROM)}.
However, the two-dimensional I-ROM~\eqref{eqn:schematic-3} has a big problem:
It is {\it not closed} since the equations for $a_1$ and $a_2$ depend on $a_3$. 
{\it This is the ROM closure problem}.

So how do we solve the ROM closure problem?
{\it The easiest way to solve the ROM closure problem is to simply ignore it}.
That is, we can simply ignore the $a_3$ contribution to the dynamics in~\eqref{eqn:schematic-3}:
\begin{equation}
\begin{aligned}
    \begin{bmatrix} 
        \dot{a}_1\\ 
        \dot{a}_2 
    \end{bmatrix} 
    = \begin{bmatrix}              
        F_1(a_1,a_2,0)\\ 
        F_2(a_1,a_2,0) 
    \end{bmatrix}.
\end{aligned} 
\label{eqn:schematic-4} 
\end{equation}
The ROM in~\eqref{eqn:schematic-4} is two-dimensional and closed (since the equations depend only on $a_1$ and $a_2$).
In Fig.~\ref{fig:schematic}, the time evolution of the solution of this two-dimensional ROM~\eqref{eqn:schematic-4} is represented as the red curve.
Of course, since in~\eqref{eqn:schematic-4} we simply ignored the $a_3$ contribution to the correct dynamics of $a_1, a_2$ given by~\eqref{eqn:schematic-3}, we incur an error, which is generally called the {\it closure error} (the red dashed lines in Fig.~\ref{fig:schematic}).

\begin{remark}[Galerkin Closure is a Relative Concept]
We note that if we start with just two ROM basis functions $\bphi_1$ and $\bphi_2$, the Galerkin ROM framework (which is presented in Section~\ref{sec:g-rom} and outlined in Algorithm~\ref{alg:g-rom}) yields a two-dimensional G-ROM that satisfies exactly the equations in~\eqref{eqn:schematic-4}.
Thus, {\it the ROM closure concept is relative to the ROM space used in the Galerkin framework}:
\begin{itemize}
    \item If we start with two basis functions, the Galerkin method yields the two-dimensional G-ROM~\eqref{eqn:schematic-4}, which is closed.
    \item If, however, we start with the larger (three-dimensional) ROM space spanned by $\bphi_1, \bphi_2$, and $\bphi_3$, the discussion in this section shows that the most accurate two-dimensional ROM obtained by a direct truncation of the three-dimensional G-ROM~\eqref{eqn:schematic-2} (i.e., the I-ROM~\eqref{eqn:schematic-3}) is not closed. 
\end{itemize}
\end{remark}

\begin{remark}[Galerkin Closure is a General Concept]
We emphasize that, although our discussion focuses exclusively on ROMs, {\it the Galerkin closure is a general concept that is associated with the classical Galerkin framework}.
Thus, there is no surprise that, over half a century, closure has been addressed in different contexts: 
large eddy simulation (LES)~\cite{BIL05}, 
variational multiscale (VMS) methods~\cite{hughes1995multiscale},
subgrid-scale (SGS) methods~\cite{guermond1999stabilization,layton2002connection},
and nonlinear Galerkin (NG) methods~\cite{foias2001navier}.
\end{remark}

\bigskip

At this point, it is probably a good idea to summarize our discussion.
As illustrated in the schematic in Fig.~\ref{fig:schematic}, the reader interested in constructing the most accurate two-dimensional G-ROM  has reached a {\it crossroads}:
\begin{itemize}
    \item On the one hand, the I-ROM~\eqref{eqn:schematic-3} is the most accurate two-dimensional ROM that we can get by using the Galerkin framework, but it is not closed.
    \item On the other hand, the G-ROM~\eqref{eqn:schematic-4} is closed, but we are incurring the closure error.
\end{itemize}
{\it This is as far as the classical Galerkin framework can take us.}
We're stuck.
So what do we do next?

The answer, as many times in numerical methods, is to take a middle of the road approach.
Specifically, we construct a {\it ROM closure model} and add it to the G-ROM~\eqref{eqn:schematic-4}:
\begin{equation}
\begin{aligned}
    \begin{bmatrix} 
        \dot{a}_1\\ 
        \dot{a}_2 
    \end{bmatrix} 
    = \begin{bmatrix}              
        F_1(a_1,a_2,0)
        + \tau_1(a_1,a_2)\\ 
        F_2(a_1,a_2,0) 
        + \tau_2(a_1,a_2) 
    \end{bmatrix},
\end{aligned} 
\label{eqn:schematic-5} 
\end{equation}
where $\tau_1(a_1,a_2), \tau_2(a_1,a_2)$ are the components of the ROM closure model, i.e., correction terms that aim at steering the inaccurate G-ROM~\eqref{eqn:schematic-4} as close as possible to the accurate (but not closed) I-ROM~\eqref{eqn:schematic-3}.
In Fig.~\ref{fig:schematic}, the time evolution of the solution of the closed ROM~\eqref{eqn:schematic-5} is represented as the green curve.

How do we construct the ROM closure model in~\eqref{eqn:schematic-5}?
We answer this question in Section~\ref{sec:d2-vms-rom}.
But first, in Section~\ref{sec:g-rom}, we present the main steps in the G-ROM construction.

\section{Galerkin ROM (G-ROM)}
	\label{sec:g-rom}

Over the past four decades, projection ROMs have been used in the numerical simulation of fluid flows~\cite{brunton2019data,hesthaven2015certified,HLB96,noack2011reduced,quarteroni2015reduced,taira2019modal}. 
In this tutorial, we exclusively consider projection ROMs that use numerical or experimental data to find the ``best" basis, which is then used together with the Galerkin method to construct the ROM.
In this section, we present the main steps in the construction of the Galerkin ROM.

To illustrate the Galerkin ROM construction, we start with a generic PDE for the dynamics of a variable of interest, $\bu$: 
\begin{eqnarray}
	{\bu}_{t}
	= 
	\bff(\bu) \, ,
	\label{eqn:pde-strong}
\end{eqnarray}
equipped with appropriate boundary conditions and initial conditions.
In Algorithm~\ref{alg:g-rom}, we list the main steps in the Galerkin ROM construction.

\begin{algorithm}[H]
	\caption{Galerkin ROM (G-ROM) Algorithm}
	\label{alg:g-rom}
	\begin{algorithmic}[1]
		\STATE{
				Use numerical or experimental data to construct modes $\{ \bphi_1, \ldots, \bphi_R \}$, which represent the recurrent spatial structures in the system~\eqref{eqn:pde-strong}.
					}
		\STATE{
				Choose the dominant modes $\{ \bphi_1, \ldots, \bphi_r \}$, $r \leq R$, as  ROM basis functions.
					}
		\STATE{
				Use a Galerkin expansion $\bur(\bx,t) = \sum_{j=1}^{r} a_j(t) \, \bphi_j(\bx)$.
					}
		\STATE{
		        Replace $\bu$ with $\bur$ in~\eqref{eqn:pde-strong}, and then on both sides of \eqref{eqn:pde-strong} take the inner product with each mode $\bphi_i, i = 1, \ldots, r$. 
		        That is, perform a Galerkin projection of the PDE~\eqref{eqn:pde-strong} onto the ROM space $\bXr:= \text{span} \{ \bphi_1, \ldots, \bphi_r \}$.  
		        The obtained {\it Galerkin ROM (G-ROM)} is of the form
                    \begin{eqnarray}
	                \overset{\bullet}{\ba}
	                =
	                \bFF(\ba),
	                \label{eqn:g-rom_alg1} 
                    \end{eqnarray}
                where $\ba(t) = (a_{i}(t))_{i=1, \ldots,r}$ is the vector of coefficients in the Galerkin expansion in step 3 
                and $\bFF$ comprises the ROM operators. 
		            }
		\STATE{
		        In the offline stage, compute the ROM operators (e.g., vectors, matrices, and tensors), which are preassembled from the ROM basis.
		            }
		\STATE{
		        In the online stage, repeatedly use the G-ROM~\eqref{eqn:g-rom_alg1} for 
		        longer time intervals. 
		            }		            
\end{algorithmic}
\end{algorithm}

\begin{remark}[ROM=d2G]
The main steps in the G-ROM~\eqref{eqn:g-rom_alg1} construction presented in Algorithm~\ref{alg:g-rom} are straightforward.
In principle, they are the same steps as those used to construct classical Galerkin methods, e.g., the finite element method (FEM).
The fundamental difference between the G-ROM and the FEM is that the former uses a {\it data-driven basis}, whereas the latter uses a universal basis (i.e., piecewise polynomials).
Thus, one could think of the projection ROMs that we discuss in this tutorial as \underline{\it data-driven Galerkin (d2G) methods}.
\label{remark:d2G}
\end{remark}

Next, we explain some of the steps in Algorithm~\ref{alg:g-rom}.

\paragraph{ROM basis (Step 1)}
To construct the ROM basis, we first collect {\it snapshots} from the simulation of the FOM.
If we are interested in time prediction (as in the numerical illustration in Section~\ref{sec:numerical-results}), the snapshots can be FEM approximations of~\eqref{eqn:pde-strong} at the time instances $t_{1}, \ldots, t_{M}$, i.e., $\bu_{h}^{1}, \ldots, \bu_{h}^{M}$, respectively. 
(If~\eqref{eqn:pde-strong} depends on parameters, we can also build a ROM basis for parameter prediction~\cite{hesthaven2015certified,quarteroni2015reduced}.)
Next, we use these snapshots to construct the modes $\{ \bphi_1, \ldots, \bphi_R \}$, which represent the recurrent spatial structures in the system described by~\eqref{eqn:pde-strong}.
Different approaches can be used to construct the ROM basis functions, e.g., 
(i) the proper orthogonal decomposition (POD)~\cite{brunton2019data,HLB96,KV01,volkwein2013proper,taira2019modal}; 
(ii) the reduced basis method (RBM)~\cite{hesthaven2015certified,quarteroni2015reduced};
(iii) the proper generalized decomposition (PGD)~\cite{chinesta2011short};
and
(iv) clustering~\cite{burkardt2006pod}.
In this tutorial, to fix ideas, we exclusively use the POD to generate the ROM basis.

For a careful presentation of the POD basis, the reader is referred to, e.g.,~\cite{HLB96} (for a physical presentation) and to~\cite{volkwein2013proper} (for a mathematical presentation).
In this paper, however, we only briefly discuss the {\it qualitative} properties of the POD basis functions, which we will later use in our numerical illustration in Section~\ref{sec:numerical-results}.
The reason for our brief qualitative discussion of the POD basis is that {\it ROM closure modeling does not depend on the particular type of ROM basis functions used}.
That is, our presentation of ROM closure modeling remains the same for any type of ROM basis used in a Galerkin framework, whether it is POD, RBM, or PGD.

The main principle used to construct the G-ROM basis can be stated as follows:
{\it Use the available snapshots to find the ROM basis that ``best'' represents the system's dynamics.}
Since this is the ``best'' basis, for certain problems, one can hope to use much fewer basis functions to construct the G-ROM than to construct, e.g., FEM models.
For example, instead of using millions or even billions of basis functions as in FEM simulations, one can hope to use tens or hundreds basis functions in the G-ROM construction.
This choice of ``best'' basis yields computational models (i.e., ROMs) whose dimension can be orders of magnitude lower than the dimension of FEM models.
(This also explains the term ``reduced'' in the ROM terminology.)

Of course, a natural question is what the ``best'' ROM basis means.
In fact, there are many proposals for the ``best'' ROM basis, and each proposal yields a different class of ROMs (e.g., POD, RBM, or PGD, to name just a few).
For example, given a set of snapshots, the POD basis is the orthonormal basis that yields the minimum projection error with respect to a chosen norm (e.g., the $L^2$ norm)~\cite{volkwein2013proper}.

\begin{figure}[h!]
	\sidecaption
 		\subfloat[ROM basis functions, $\varphi_1$ and $\varphi_{10}$]{\includegraphics[width=.9\textwidth]{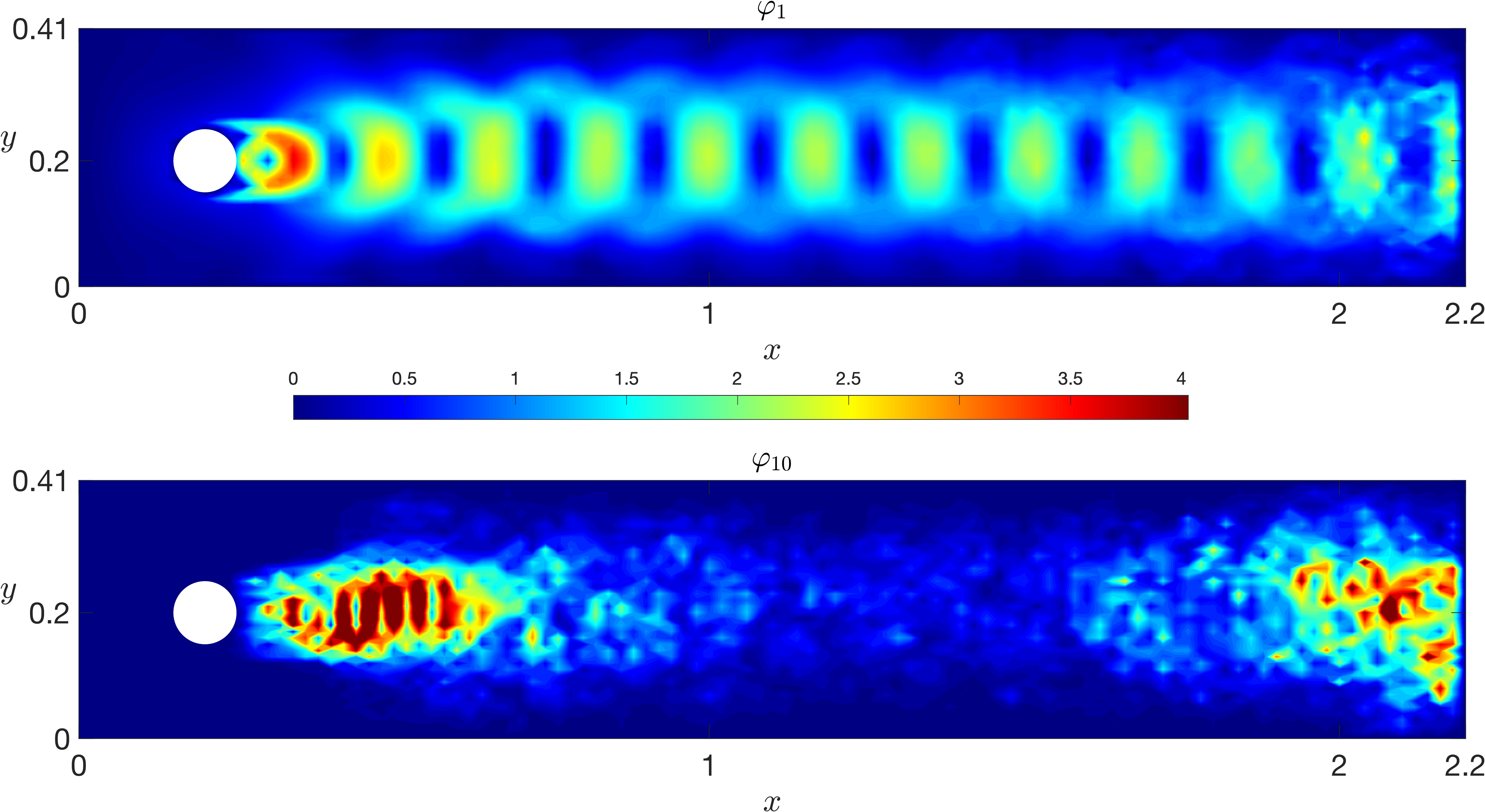}}
 		\\
  		\subfloat[FEM basis functions, $\phi^h_1$ and $\phi^h_{10}$]{
  		\includegraphics[width=.9\textwidth]{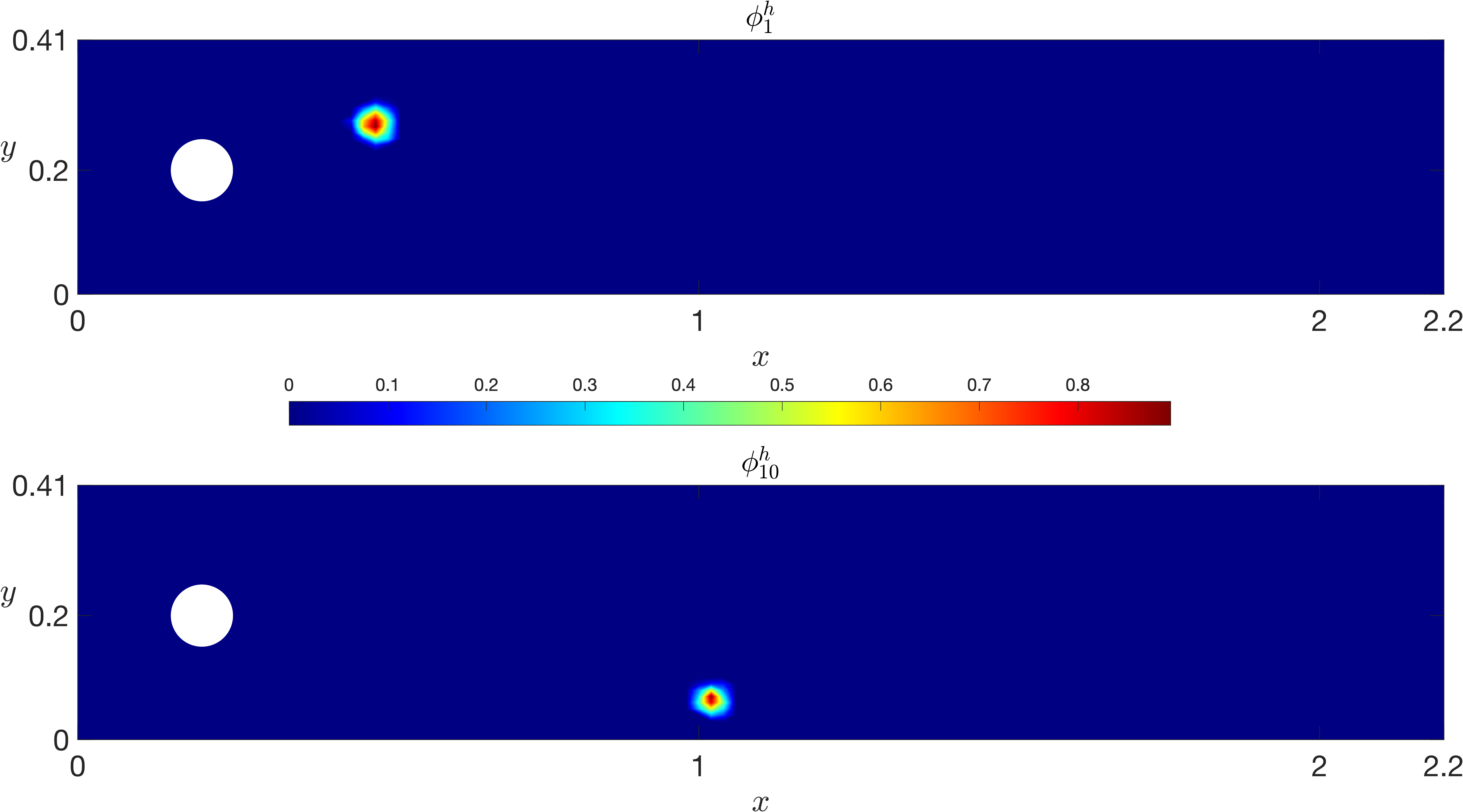}
  		}
	\caption{
	    2D flow past a circular cylinder:
		(a) ROM basis functions $\bphi_{1}$ and $\bphi_{10}$.
		(b) FEM basis functions $\phi_{1}^h$ and $\phi_{10}^h$.
		Note that the ROM basis functions are fundamentally different from the FEM basis functions:
		The former have global support, whereas the latter have local support.
	}
	\label{fig:basis-nse}
\end{figure}

\begin{figure}[h!]
	\sidecaption
  		\subfloat{
  		\includegraphics[width=.5\textwidth]{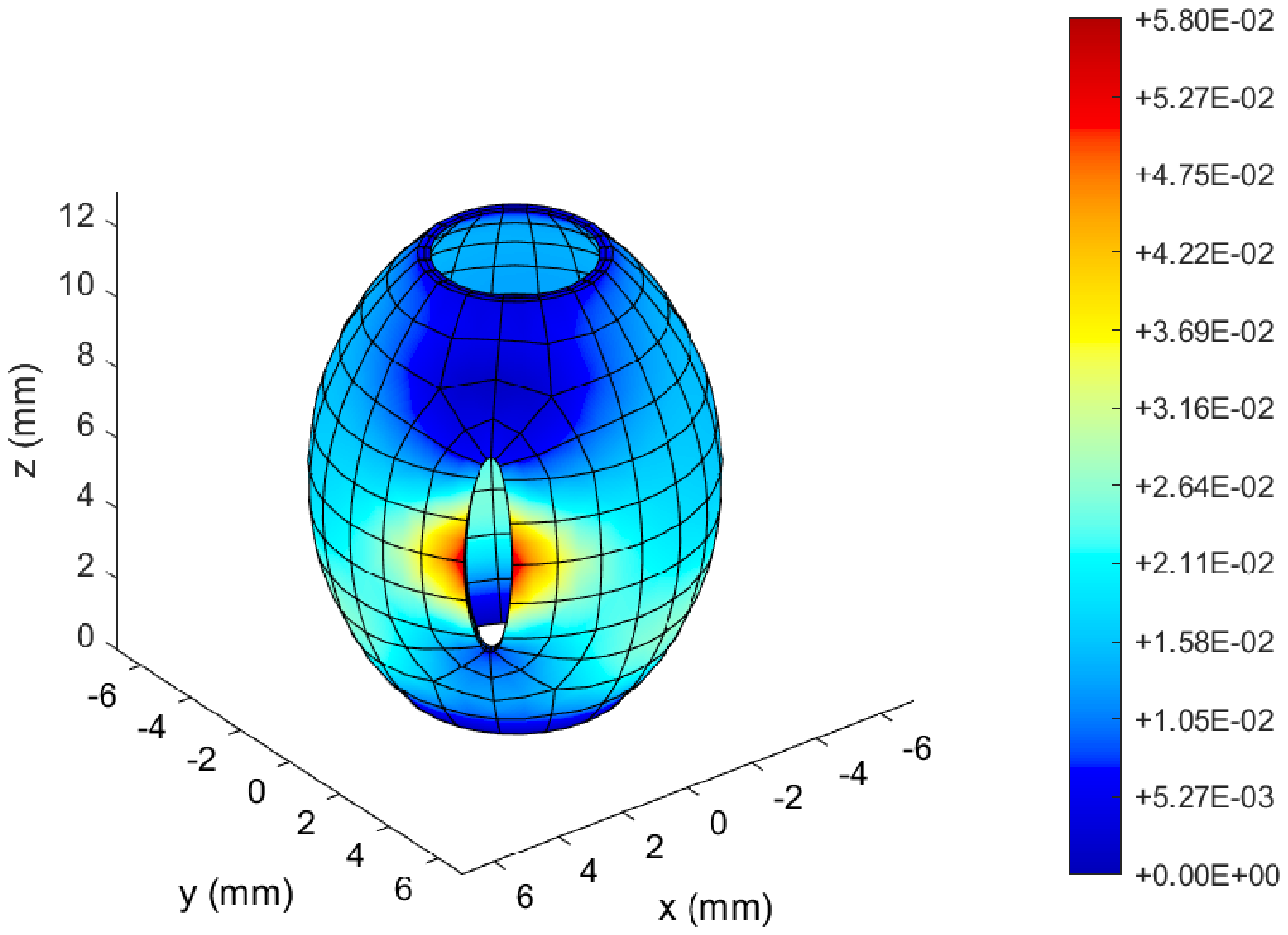}
  		}
  		\hfill
  		\subfloat{
  		\includegraphics[width=.5\textwidth]{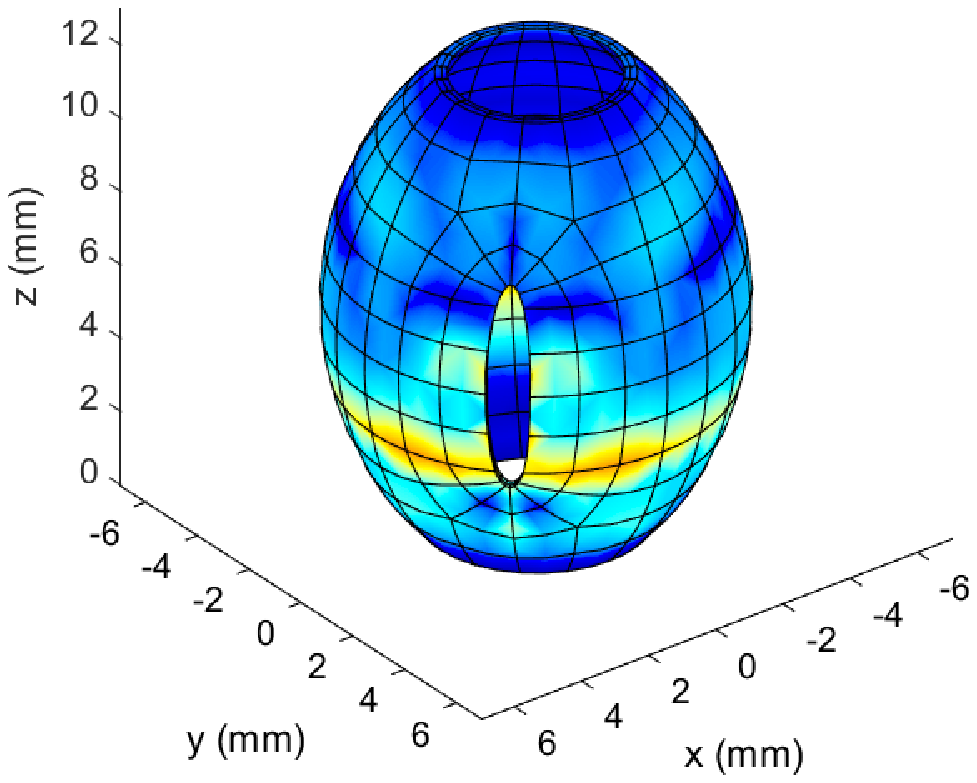}
  		}
	\caption{
	    Soft tissue modeling:
		ROM basis functions $\bphi_{1}$ and $\bphi_{10}$.
	}
	\label{fig:basis-tissue}
\end{figure}

However, independent of the approach used to construct them, the ROM basis functions generally share several qualitative features.
To illustrate this, in Fig.~\ref{fig:basis-nse} we plot two POD basis functions, $\bphi_{1}$ and $\bphi_{10}$, and two FEM basis functions, $\boldsymbol{\phi}_{1}^{h}$ and $\boldsymbol{\phi}_{10}^{h}$, for a 2D flow past a circular cylinder~\cite{mou2021data}. 
One can clearly see the significant differences between the POD basis functions (top two plots) and the FEM basis functions (bottom two plots).
Indeed, the POD basis functions have global support (i.e., they can be nonzero over the entire computational domain), whereas the FEM basis functions have local support (i.e., they are one at one mesh point and zero everywhere else).
To further illustrate the different characteristics of the POD basis, in Fig.~\ref{fig:basis-tissue} we plot two POD basis functions, $\bphi_{1}$ and $\bphi_{10}$, for soft tissue modeling~\cite{snyder2022data}.
Comparing these two POD basis functions with the POD basis functions in the top two plots of Fig.~\ref{fig:basis-nse}, we can clearly see that different physical systems (i.e., the soft tissue in Fig.~\ref{fig:basis-tissue} and the flow in Fig.~\ref{fig:basis-nse}) yield fundamentally different POD basis functions.
We emphasize that this is in complete contrast with classical numerical methods, such as the FEM.
Indeed, the FEM basis functions are {\it universal} basis functions, i.e., they have the same shape (piecewise polynomials and local support) for all the problems.
In contrast, the POD basis functions (and ROM basis functions in general) change their shape when we change the problem.
This can be clearly seen by comparing the top two plots of Fig.~\ref{fig:basis-nse} with the plots of Fig.~\ref{fig:basis-tissue}. 

\vspace*{1.0cm} 
\paragraph{Galerkin ROM construction (Steps 2-6)}

To illustrate the G-ROM construction, we use the Navier-Stokes equations (NSE) as a mathematical model:
\vspace*{-0.1cm}
 \begin{eqnarray}
     && \frac{\partial \bu}{\partial t}
     - Re^{-1} \Delta \bu
     + \bu \cdot \nabla \bu
     + \nabla p
     = {\bf 0} \, ,
     \label{eqn:nse-1} \\
     && \nabla \cdot \bu
     = 0 \, ,
     \label{eqn:nse-2}
	\\[-0.6cm]
	\nonumber
 \end{eqnarray}
 where $\bu$ is the velocity, $p$ the pressure, and $Re$ the Reynolds number. 
We consider the NSE posed on a bounded spatial domain in either $\mathbb{R}^2$ or $\mathbb{R}^3$, and supplemented with homogeneous Dirichlet boundary conditions and an appropriate initial condition.
The NSE~\eqref{eqn:nse-1}--\eqref{eqn:nse-2} can be cast in the general form~\eqref{eqn:pde-strong} by choosing $\bff(\bu) = Re^{-1} \Delta \bu - \bu \cdot \nabla \bu$ (after applying the Leray projection, which maps the vector field into the divergence-free subspace of the underlying state space) \cite{temam2001navier}. 

To construct the G-ROM for the NSE, we follow Steps 2-6 in Algorithm~\ref{alg:g-rom}.
That is, we choose the first $r$ basis functions from the modes constructed in Step 1, use a Galerkin truncation
\begin{eqnarray}
    \bur(\bx,t) 
    = \sum_{j=1}^{r} a_j(t) \, \bphi_j(\bx) ,
    \label{eqn:galerkin-expansion}
\end{eqnarray}
replace $\bu$ with $\bur$ in the NSE~\eqref{eqn:nse-1}, and  project the resulting PDE onto the ROM space, $\bX^r$. 
Furthermore, we apply the divergence theorem to the diffusion term and the pressure term.
This yields the G-ROM~\cite{mou2021data}:
\vspace*{-0.1cm}
\begin{eqnarray}
	\overset{\bullet}{\ba} 
	= A \, \ba 
	+ \ba^{\top} \, B \, \ba ,
\label{eqn:g-rom}
	\\[-0.6cm]
	\nonumber
\end{eqnarray}
where 
$\ba(t)$ is the vector of unknown coefficients $a_j(t), 1 \leq j \leq r$ in the Galerkin expansion~\eqref{eqn:galerkin-expansion}.
The ROM operator $A$ in~\eqref{eqn:g-rom} is an $r \times r$ matrix that corresponds to the diffusion term in the NSE (i.e., $- Re^{-1} \Delta \bu$) and has entries 
\begin{eqnarray}
 	A_{im}
	= - Re^{-1} \, \left( \nabla \bphi_m , \nabla \bphi_i \right) ,
	\quad 1 \leq i, m \leq r \, ,
    \label{eqn:rom-operator-A}   
\end{eqnarray}
where $(\cdot , \cdot)$ denotes the $L^2$ inner product.
The ROM operator $B$ in~\eqref{eqn:g-rom} is an $r \times r \times r$ tensor that corresponds to the nonlinear term in the NSE (i.e., $\bu \cdot \nabla \bu$) and has entries
\begin{eqnarray}
	B_{imn}
	= - \bigl( \bphi_m \cdot \nabla \bphi_n , \bphi_i \bigr) \, ,
	\ 
	\quad 1 \leq i, m, n \leq r \, .
    \label{eqn:rom-operator-B}   
\end{eqnarray}

We note that the pressure term in the G-ROM~\eqref{eqn:g-rom} vanishes since we assumed that the ROM modes are discretely divergence-free (which is the case if, e.g., the snapshots are discretely divergence-free).
ROMs that provide a pressure approximation are discussed in, e.g.,~\cite{decaria2020artificial,hesthaven2015certified}.

Once the matrix $A$ and tensor $B$  are assembled in the offline stage, the G-ROM~\eqref{eqn:g-rom} is a relatively low-dimensional, efficient dynamical system that can be used in the online stage for longer time intervals (or more parameter values, e.g., $Re$~\cite{hesthaven2015certified,quarteroni2015reduced}).


\section{The Closure Problem and Its Solution: The Closure Model}
	\label{sec:closure-problem-closure-model}

This section has two goals:
In Section~\ref{sec:closure-problem}, we motivate the need for ROM closure modeling in the under-resolved regime, i.e., we describe the ROM closure problem.
In Section~\ref{sec:closure-model}, we show how to solve the ROM closure problem, i.e., we show how to construct a ROM closure model.
To this end, we give the definition of the ROM closure model, show that using the exact closure model (i.e., using the ideal ROM) increases the ROM accuracy, and finally outline the main steps in the ROM closure model construction.

\subsection{The Closure Problem}
	\label{sec:closure-problem}

The G-ROM~\eqref{eqn:g-rom} constructed in Section~\ref{sec:g-rom} is appealing from the computational point of view: 
The G-ROM can significantly reduce the dimension (and, thus, the computational cost) of classical numerical discretization (e.g., FEM) models by orders of magnitude.
So one can ask the following natural question:

\begin{question}{Q1}
	What is wrong with G-ROM?
\end{question}

The short answer to Q1 is: 
It depends on the resolution.
Specifically:

\begin{important}{A1}
	It depends on whether we are in the \textit{resolved} regime or the \textit{under-resolved} regime.
\end{important}

\begin{itemize} \itemsep10pt
	\item In the \textit{resolved} regime (i.e., when there are enough ROM basis functions $\{ \bphi_{1}, \ldots, \bphi_{r} \}$ to accurately represent the underlying dynamics), the G-ROM produces accurate results. 
	\item In the \textit{under-resolved} regime (i.e., when there are not enough ROM basis functions $\{ \bphi_{1}, \ldots, \bphi_{r} \}$ to accurately represent the underlying dynamics), the G-ROM produces inaccurate results. 
\end{itemize}

But then one can ask the following questions:

\begin{question}{Q2}
	Why is the under-resolved regime important?
	Why do we need to worry about it?
\end{question}

\begin{important}{A2}
	Many important applications (e.g., atmospheric boundary layer flows, digital twins of wind farms, and anisotropic and heterogeneous biological tissues) are centered around {\it multiscale} systems that require a large number of ROM basis functions.
	However, to ensure a low computational cost in these applications, under-resolved G-ROMs are generally used.
\end{important}

\subsection{The Closure Model}
	\label{sec:closure-model}

In Section~\ref{sec:closure-problem}, we defined the ROM closure problem and we explained why it is important.
In this section, we present the solution to the ROM closure problem.
That is, we answer the following question:
\begin{question}{Q3}
	What is the solution to the closure problem?
\end{question}

\begin{important}{A3}
	The solution to the closure problem is the closure model.
	That is, replace the G-ROM~\eqref{eqn:g-rom} with 
	\begin{equation} 
		\overset{\bullet}{\ba} 
		= \bFF(\ba) 
		+ \btau(\ba), 
		\label{eqn:g-rom-tau}
	\end{equation}	
where $\btau(\ba)$ is the closure model, which represents the effect of the discarded ROM modes $\{ \bphi_{r+1}, \ldots, \bphi_{R} \}$ on the ROM dynamics.
\end{important}

Note that A3 is a vague definition, which begs the following questions:
What exactly does ``model the effect'' mean?
What exactly does $\btau(\ba)$ in~\eqref{eqn:g-rom-tau} actually model?

Answering these natural questions is not straightforward.
To do so, we need to \underline{\textit{extend the Galerkin framework}}.
This sounds like a daunting task, but it turns out to be relatively simple.
The ``trick'' is to \textit{rethink the space} we use in the Galerkin framework:

In the resolved regime, the ROM space $\bX^{r} := \text{span}\{ \bphi_{1}, \ldots, \bphi_{r} \}$ is the only space we will ever need, since everything happens in $\bX^{r}$.
Thus, in the resolved regime, G-ROM should (and generally does) work just fine.

However, in the under-resolved regime we need {\it two} spaces:
(i) the {\it resolved space} $\bX^{r}$, and (ii) the {\it unresolved space} $\bX^{r'} := \text{span}\{ \bphi_{r+1}, \ldots, \bphi_{R} \}$.
To keep the ROM dimension (and, therefore, its computational cost) low, we want to work in the resolved space, $\bX^{r}$.
However, to increase the ROM accuracy, we should do our best to model the contribution to the ROM dynamics made  by the dynamics in the unresolved space, $\bX^{r'}$. But this sounds like a lot of work (both in terms of modeling and computation).
So the following is a natural question:

\begin{question}{Q4}
	Does $\bX^{r'}$ have a significant effect on the ROM dynamics?
\end{question}

\begin{important}{A4}
	Yes.
\end{important}

The answer A4 is simple.
In section~\ref{sec:i-rom}, we introduce the ideal ROM, which adds the {\it exact} closure term to the classical G-ROM.
The ideal ROM results clearly show why the effect of $\bX^{r'}$ should be modeled.
Specifically, we show that the ideal ROM results are dramatically more accurate than the G-ROM results.
Thus, we conclude that modeling the exact ROM closure term is beneficial to ROM accuracy.

\subsubsection{The Ideal ROM (I-ROM)}
	\label{sec:i-rom}

To present the ideal ROM, we first need to define the spaces of resolved ROM scales (i.e., $\bX^{r}$) and unresolved ROM scales (i.e., $\bX^{r'}$).
To this end, we extend the variational multiscale (VMS) framework proposed by Hughes and his group two decades ago in the FEM context.
We note, however, that there are other ways of defining the spaces of resolved and unresolved ROM scales, e.g., spatial filtering~\cite{mou2021data}.

First, we leverage the orthonormality of the ROM basis functions and construct the two orthogonal spaces, $\bX^{r}$ and $\bX^{r'}$, as follows:
\begin{eqnarray}
    \bX^{r} 
    := \text{span}  \{ \bphi_1, \ldots, \bphi_{r} \}
    \qquad \text{and} \qquad 
    \bX^{r'} 
    := \text{span}  \{ \bphi_{r+1}, \ldots, \bphi_{R} \} .
    \label{eqn:resolved-unresolved-decomposition}
\end{eqnarray}
The space $\bX^{r}$ represents the space of the {\it resolved ROM scales}, i.e., the spatial scales that are explicitly approximated by a given  $r$-dimensional ROM.
In contrast, the space $\bX^{r'}$ represents the space of the {\it unresolved ROM scales}, i.e., the spatial scales that are not explicitly approximated by the chosen ROM.
We note that since the ROM basis functions are generally ordered from the most important to the least important (with respect to a physical criterion, e.g., kinetic energy~\cite{HLB96}), the decomposition in~\eqref{eqn:resolved-unresolved-decomposition} is natural.
We also note that since we are concerned with the under-resolved regime that often occurs in practical applications, we consider the case when $r \ll R$.

The next step in the construction of the ideal ROM is to extend the Galerkin framework to the space $\bX^{R} := \bX^{r} \oplus \bX^{r'}$, which is the {\it maximal ROM space} (i.e., the space spanned by all the snapshots). 
Thus, we use the ROM approximation of both resolved and unresolved scales, i.e., we utilize  $\bu_{R} \in \bX^{R}$  defined as
\begin{eqnarray}
	\bu_{R}
	= \sum_{j=1}^{R} a_{j} \, \bphi_{j}
	= \sum_{j=1}^{r} a_{j} \, \bphi_{j}
	+ \sum_{j=r+1}^{R} a_{j} \, \bphi_{j} = \bu_{r}
	+ \bu' \, ,
	\label{eqn:u_R}
\end{eqnarray}
where $\bu_{r} \in \bX^{r}$ represents the {\it resolved ROM component} of $\bu$, and $\bu' \in \bX^{r'}$ represents the {\it unresolved ROM component} of $\bu$.
Next, we plug $\bu_R$ in the generic equation~\eqref{eqn:pde-strong}, project the resulting equation onto $\bX^{r}$, and use the ROM basis orthogonality  to show that  
$\bigl( \bu_{R,t} \, , \bphi_{i} \bigr) = \bigl( \bu_{r,t} \, , \bphi_{i} \bigr), \  \forall \, i = 1, \ldots, r$, where $\bu_{R,t}$ and $\bu_{r,t}$ are the time derivatives of $\bu_{R}$ and $\bu_{r}$, respectively. 
Following these steps, we obtain the {\it ideal ROM (I-ROM)}:
\begin{eqnarray}
		\bigl( \bu_{r,t} \, , \bphi_{i} \bigr) 
		= 
		\bigl( \bff(\bu_{r}) \, , \bphi_{i} \bigr) 
		+
		\underbrace{
			\bigl( \bff(\bu_{R}) \, , \bphi_{i} \bigr)
			-
			\bigl( \bff(\bu_{r}) \, , \bphi_{i} \bigr)
		}_{\btau^{I-ROM} = \text{ ideal ROM closure term}}
		\ , 
		\  
		\forall \, i = 1, \ldots, r .
	\label{eqn:i-rom}
\end{eqnarray}
The last two terms in~\eqref{eqn:i-rom} yield the {\it ideal ROM closure term}, $\btau^{I-ROM}$, which represents the effect of the discarded ROM modes $\{ \bphi_{r+1}, \ldots, \bphi_{R} \}$ onto the dynamics of the resolved ROM scales, $\bur$.
Using the expansion~\eqref{eqn:u_R}, the I-ROM~\eqref{eqn:i-rom} can be written as the following dynamical system for the vector of ROM coefficients of the resolved scales: 	
\begin{equation} 
	\overset{\bullet}{\ba} 
	= \bFF(\ba) 
	+ \btau^{I-ROM}(a_{1}, \ldots, a_{r},a_{r+1}, \ldots, a_{R}). 
	\label{eqn:i-rom-ds}
\end{equation}	

The above discussion clearly shows that, from a mathematical point of view, the correct equations satisfied by the coefficients of the resolved ROM scales are the I-ROM equations~\eqref{eqn:i-rom-ds} instead of the G-ROM equations~\eqref{eqn:g-rom}. 
However, we need to ask ourselves whether this mathematical framework has a {\it practical impact} (i.e., we need to ask question Q4).
Specifically, we need to check whether the I-ROM results are better than the G-ROM results.

\begin{figure}[h!]
	\sidecaption
	\centering
    \includegraphics[width=.9\textwidth]{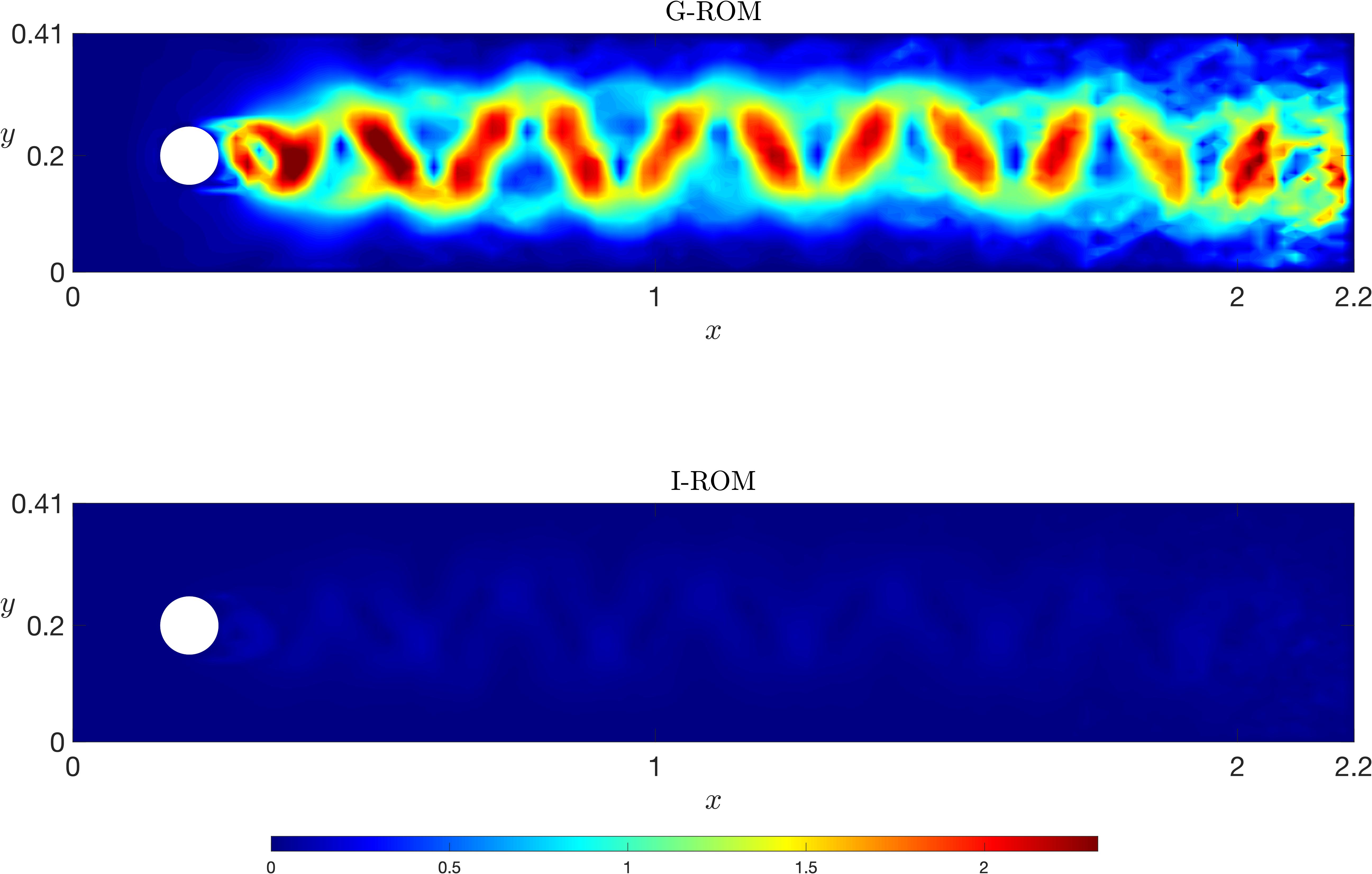}
	\caption{
        2D flow past a circular cylinder. $L^2$ norm of the error, $\| \bu^{FOM}-\bu^{ROM} \|_{L^2}$, for G-ROM~\eqref{eqn:g-rom} (top) and  I-ROM~\eqref{eqn:i-rom} (bottom).
        The I-ROM error is significantly lower than the G-ROM error, which illustrates the potential benefit of ROM closure modeling.
	}
	\label{fig:i-rom-nse}
\end{figure}

In Fig.~\ref{fig:i-rom-nse}, we present results for the I-ROM~\eqref{eqn:i-rom-ds} and the G-ROM~\eqref{eqn:g-rom} in the numerical simulation of a two-dimensional flow past a circular cylinder.
These plots clearly show that the I-ROM performs significantly better than the classical G-ROM.
Thus, these results suggest that including a model for the I-ROM closure term, $\btau^{I-ROM}$, could increase the ROM accuracy.

\begin{remark}[The Closure Model Increases Accuracy]
There is a lot of confusion in the ROM community (and not only) regarding the role of the closure model.
In this section, we tried to emphasize that the main role of the ROM closure model is to increase the {\it accuracy} of the G-ROM.
Indeed, in equation~\eqref{eqn:g-rom-tau}, adding the closure term, $\btau(\ba)$, to the classical G-ROM yields a more accurate model (in the extended Galerkin framework).

That being said, in many important practical applications (e.g., convection-dominated flows), the G-ROM's inaccuracy often manifests itself in the form of {\it spurious numerical oscillations}.
Thus, a popular misconception (at least in computational fluid dynamics) is that the only role of the ROM closure model is to eliminate/alleviate these numerical oscillations, i.e., to increase the numerical stability of the G-ROM.

However, we emphasize that, while numerical stability of the model is {\it necessary} (indeed, if the model is accurate, then it has to be stable), it is {\it not sufficient}.
For example, we can add a very large stabilization term to the classical G-ROM.
This, most likely, will stabilize the model, but will also degrade its accuracy.

To summarize, we emphasize that ROM closure modeling is not simply about adding numerical stabilization. 
Instead, ROM closure modeling is about adding the {\it ``right''} amount of numerical stabilization (i.e., the amount of stabilization that makes the model accurate).
\end{remark}

\subsubsection{Closure Model Construction}
    \label{sec:closure-model-construction}

The I-ROM results in Section~\ref{sec:i-rom} clearly show that the effect of $\bX^{r'}$ should be modeled.
We emphasize, however, that the I-ROM itself does {\it not} represent a practical solution since it depends on the coefficients of the discarded ROM modes, $a_{r+1}, \ldots, a_{R}$, which we do not model in our ROM (since we work in $\bX^{r}$).

\begin{question}{Q5}
	How do we make the I-ROM~\eqref{eqn:i-rom-ds} practical?
\end{question}

\begin{important}{A5}
	We construct a closure model, $\btau$,  which is an approximation in $\bX^{r}$ of the I-ROM closure term, $\btau^{I-ROM}$:
	\begin{equation}
		\btau^{I-ROM}(a_{1}, \ldots, a_{r}, a_{r+1}, \ldots, a_{R})
		\approx \btau(a_{1}, \ldots, a_{r}).
		\label{eqn:i-rom-closure}
	\end{equation}
\end{important}

Since $\btau$ in~\eqref{eqn:i-rom-closure} lives in $\bX^{r}$, it can be computed with the available ROM data, and, thus, can be used in practical computations.

\begin{remark}[Closure=Correction]
Equation~\eqref{eqn:g-rom-tau} shows that the closure model, $\btau$, in~\eqref{eqn:i-rom-closure} can be interpreted as a correction term that is added to the G-ROM~\eqref{eqn:g-rom} to correct its dynamics  in $\bX^{R}$.
So do we really need I-ROM in order to construct the closure model?
In Section~\ref{sec:d2-vms-rom}, we will show that the I-ROM is needed when we construct data-driven ROM closures.
Furthermore, we note that the I-ROM derivation explains the closure model terminology.
Indeed, $\btau^{I-ROM}(a_{1}, \ldots, a_{r},a_{r+1}, \ldots, a_{R})$ shows that the I-ROM~\eqref{eqn:i-rom} is closed in $\bX^{R}$, but not in $\bX^{r}$. 
\end{remark}

ROM closure models are of three types:
(i) {\it Functional}, which use physical insight to construct the closure model.
(ii) {Structural}, which use mathematical tools.
(iii) {\it Data-driven}, which use available data.
The three types of ROM closure models are surveyed in~\cite{ahmed2021closures}.
In this tutorial, we take a different approach and, for clarity of presentation, focus on data-driven approaches, which have experienced a tremendous development over the last few years.
Specifically, in the next section, we present the data-driven variational multiscale ROM closure model.

\section{The Data-Driven Variational Multiscale ROM (D2-VMS-ROM)}
	\label{sec:d2-vms-rom}

In this section, we illustrate how data-driven modeling can be leveraged to construct the ROM closure model.
Specifically, we outline the main steps in the construction of one data-driven ROM closure model, i.e., the data-driven variational multiscale ROM (D2-VMS-ROM) that was proposed in~\cite{mou2021data} (see also~\cite{xie2018data}).
To this end, we follow the presentation in Section 2.3 in \cite{mou2021data} to construct the two-scale D2-VMS-ROM.
(We note that a three-scale D2-VMS-ROM was also proposed and tested in \cite{mou2021data}.)

To build the D2-VMS-ROM, we start with the I-ROM~\eqref{eqn:i-rom-ds}.
As explained in answer A5, to construct the ROM closure model we need to find an approximation $\btau(a_{1}, \ldots, a_{r})$ for the I-ROM closure term in~\eqref{eqn:i-rom-ds}, $\btau^{I-ROM}(a_{1}, \ldots, a_{r}, a_{r+1}, \ldots, a_{R})$.
The construction of the data-driven ROM closure model consists of two steps: (i) postulating a model form ansatz; and (ii) solving a least squares problem to determine the coefficients of the model form.
Next, we outline these two steps.

\paragraph{Model Form Ansatz}
The first step in the construction of the data-driven ROM closure model is to {\it postulate a model form (ansatz)}.
Specifically, we approximate the I-ROM closure term $\btau^{I-ROM}$ with ${\boldsymbol{g}}(\bu_r)$, where ${\boldsymbol{g}}$ is a {\it generic} function whose coefficients/parameters still need to be determined:
\begin{eqnarray}
\boxed{
	\btau_{i}^{I-ROM} 
	\stackrel{\eqref{eqn:i-rom}}{=}  
			\bigl( \bff(\bu_{R}) \, , \bphi_{i} \bigr)
			-
			\bigl( \bff(\bu_{r}) \, , \bphi_{i} \bigr)
		\approx
		\bigl( {\boldsymbol{g}}(\bu_{r}) \, , \bphi_{i} \bigr) \,, \quad i = 1, \ldots, r.
	}
	\label{eqn:g}
\end{eqnarray}

\paragraph{Least Squares Problem}
To determine the coefficients/parameters in ${\boldsymbol{g}}$ used in~\eqref{eqn:g}, in the offline  stage, we solve the following low-dimensional {\it least squares problem}:
\begin{eqnarray}
	\boxed{
	 \begin{aligned}
	\min_{ {\boldsymbol{g}}\text{ parameters}} \ \sum_{j=1}^{M} \biggl\| 
	\bigl[ 
			\bigl( \bff(\bu_{R}^{FOM}(t_j)) \, , \bphi_{i} \bigr)
			&-
			\bigl( \bff(\bu_{r}^{FOM}(t_j)) \, , \bphi_{i} \bigr)
		\bigr] 
	 \\&-
	\bigl( {\boldsymbol{g}}(\bu_{r}^{FOM}(t_{j})) \, , \bphi_{i} \bigr) \biggr\|^2  ,
	 \end{aligned}
	}
	\label{eqn:least-squares}
\end{eqnarray}
where $\bu_{R}^{FOM}$ and $\bu_{r}^{FOM}$ are obtained from the FOM data, and $M$ is the number of snapshots.
Once ${\boldsymbol{g}}$ is determined, the I-ROM~\eqref{eqn:i-rom} with the I-ROM closure term replaced by ${\boldsymbol{g}}$ yields the {\it data-driven VMS-ROM (D2-VMS-ROM)}:
\begin{eqnarray}
	\boxed{
	\bigl( \bu_{r,t} \, , \bphi_{i} \bigr) 
	= \bigl( \bff(\bu_{r}) \, , \bphi_{i} \bigr) 
	+ \bigl( {\boldsymbol{g}}(\bu_{r}) \, , \bphi_{i} \bigr) ,
	}
	\qquad 
	i = 1, \ldots, r.
	\label{eqn:d2-vms-rom}
\end{eqnarray}

We emphasize that we have a lot of {\it flexibility} in choosing the model form ansatz~\eqref{eqn:g} in the D2-VMS-ROM.
For example, for the NSE, we can choose the following model form:
$\forall \, i = 1, \ldots, r,$
\begin{eqnarray}
    \bigl( {\boldsymbol{g}}(\bu_{r}) \, , \bphi_{i} \bigr) 
	= \bigl( \tilde{A} \, \ba + \ba^\top \tilde{B} \, \ba \bigr)_{i} \ ,
	\label{eqn:ansatz-nse}
\end{eqnarray}
where, for computational efficiency, we assume that the structures of ${\boldsymbol{g}}$  and $\bff$ are similar.
Thus, in the least squares  problem~\eqref{eqn:least-squares}, we solve for all the entries in the $r \times r$ matrix $\tilde{A}$ and the $r \times r \times r$ tensor $\tilde{B}$.

The least squares problem~\eqref{eqn:least-squares} is {\it low-dimensional} since there are only $(r^2 + r^3)$ entries in $\tilde{A}$ and $\tilde{B}$ to be optimized, 
and $r$ is small.
Thus,~\eqref{eqn:least-squares} can be efficiently solved in the offline stage.
For the NSE, the D2-VMS-ROM~\eqref{eqn:d2-vms-rom} takes the form
\begin{eqnarray}
	\boxed{
	\overset{\bullet}{\ba} 
	= ( A + \tilde{A}) \ba 
	+ \ba^\top (B + \tilde{B}) \ba \, ,
	}
	\label{eqn:d2-vms-rom-ds}
\end{eqnarray}
where $A$ and $B$ are the G-ROM  operators in~\eqref{eqn:g-rom}, and $\tilde{A}$ and $\tilde{B}$ are the VMS-ROM closure operators in~\eqref{eqn:ansatz-nse}. 

\begin{remark}[Physical Constraints]
To improve the D2-VMS-ROM accuracy, one can use physical constraints when solving the least squares problem~\eqref{eqn:least-squares} to find the entries of the VMS-ROM closure operators $\tilde{A}$ and $\tilde{B}$.
Numerical experiments have shown that imposing physical constraints can indeed increase the D2-VMS-ROM accuracy~\cite{mohebujjaman2019physically}.
\end{remark}

In Algorithm~\ref{alg:d2-rom-closure}, we list the main steps in the construction of ROMs equipped with data-driven closure models.

\begin{algorithm}[H]
	\caption{Data-Driven ROM Closure Algorithm}
	\label{alg:d2-rom-closure}
	\begin{algorithmic}[1]
		\STATE{
				Use numerical or experimental data to construct modes $\{ \bphi_1, \ldots, \bphi_R \}$, which represent the recurrent spatial structures in the system.
					}
		\STATE{
				Choose the dominant modes $\{ \bphi_1, \ldots, \bphi_r \}$, $r \leq R$, as  ROM basis functions.
					}
		\STATE{
				Use a Galerkin expansion $\bu_{R}(\bx,t) = \sum_{j=1}^{R} a_j(t) \, \bphi_j(\bx)$.
					}
		\STATE{
		        Replace $\bu$ with $\bu_{R}$ in~\eqref{eqn:pde-strong}.
		            }
		\STATE{
		        Use a Galerkin projection of the PDE obtained in step 4 onto the space of resolved ROM scales $\bXr := \text{span} \{ \bphi_1, \ldots, \bphi_r \}$ to obtain the {\it ideal ROM (I-ROM)}:
                    \begin{eqnarray}
	                \overset{\bullet}{\ba}
	                =
	                \bFF(\ba)
	                + \btau^{I-ROM}
	                \label{eqn:alg-d2-vms-rom-i-rom} 
                    \end{eqnarray}
                where $\ba(t) = (a_{i}(t))_{i=1, \ldots,r}$ is the vector of coefficients in the Galerkin expansion in step (3), $\bFF$ comprises the G-ROM operators, and $\btau^{I-ROM}$ is the ideal ROM closure term defined in~\eqref{eqn:i-rom}.
		            }
		\STATE{
		        In the offline stage: 
		            \begin{itemize}
		                \item Compute the G-ROM operators (e.g., vectors, matrices, and tensors), which are preassembled from the ROM basis.
		                \item Choose a model form $\boldsymbol{g}$ for $\btau^{I-ROM}$ in~\eqref{eqn:alg-d2-vms-rom-i-rom}.
		                \item Solve the least squares problem~\eqref{eqn:least-squares} to find the parameters in the model form.
		                \item Compute $\bG(\ba)$, 
		                which comprises the ROM closure operators corresponding to the model form $\boldsymbol{g}$ for $\btau^{I-ROM}$.
		                \item Replace the I-ROM~\eqref{eqn:alg-d2-vms-rom-i-rom} with the data-driven ROM closure model
		                \begin{eqnarray}
	                        \overset{\bullet}{\ba}
	                        =
	                        \bFF(\ba)
	                        + \bG(\ba)
	                        \label{eqn:alg-d2-vms-rom-closure} 
                    \end{eqnarray}
		            \end{itemize}
		            }
		\STATE{
		        In the online stage, repeatedly use the data-driven ROM closure~\eqref{eqn:alg-d2-vms-rom-closure} for various parameter settings and/or longer time intervals.
		            }		            
\end{algorithmic}
\end{algorithm}

\section{ROM Closures in Action: Numerical Results}
    \label{sec:numerical-results}

In the previous sections, we tried to convince the reader that ROM closures are important since they significantly increase the ROM accuracy in the under-resolved regime.
We note, however, that all our arguments have been {\it mathematical} arguments.
Thus, we can ask the following natural question:

\begin{question}{Q6}
	Do ROM closures work in practice?
\end{question}

The answer to Q6 is simple:
\begin{important}{A6}
	Yes!
\end{important}

The answer A6 is elaborated in the survey in~\cite{ahmed2021closures}, which presents a plethora of examples of under-resolved ROM simulations of complex dynamics (e.g., turbulent flows) in which ROM closures significantly increase the accuracy at a modest computational overhead.

\begin{figure}[h!]
	\sidecaption
    	    \includegraphics[width=.9\textwidth]{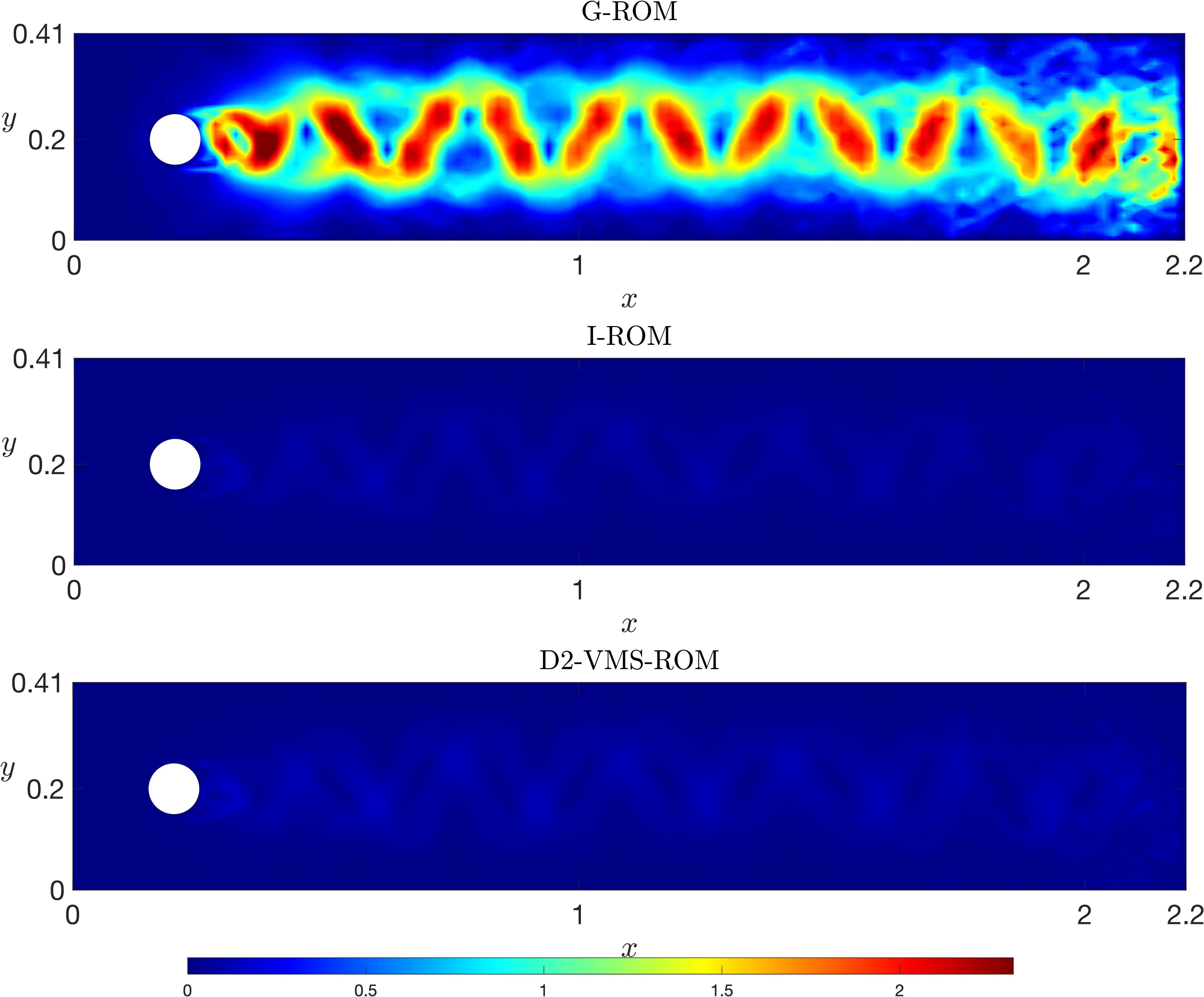}
	\caption{
        2D flow past a circular cylinder. $L^2$ norm of the error, $\| \bu^{FOM}-\bu^{ROM} \|_{L^2}$, for G-ROM~\eqref{eqn:g-rom} (top), I-ROM~\eqref{eqn:i-rom} (middle), and D2-VMS-ROM~\eqref{eqn:d2-vms-rom} (bottom).
        The D2-VMS-ROM error is significantly lower than the G-ROM error, which illustrates the benefit of ROM closure modeling.
        Also note that, in this case, the D2-VMS-ROM error almost reaches the theoretical lower bound given by the I-ROM error. 
	}
	\label{fig:d2-vms-rom-nse}
\end{figure}

In this section, for clarity of presentation, we illustrate how a specific ROM closure model (i.e., the D2-VMS-ROM outlined in Section~\ref{sec:d2-vms-rom}) increases the ROM accuracy for the 2D flow past a circular cylinder~\cite{mou2021data}, which is a simple test problem commonly used in the ROM community.  
(We note, however, that the D2-VMS-ROM was successfully used for challenging test problems, e.g., turbulent channel flow~\cite{mou2021data-phd} and the quasi-geostrophic equations~\cite{mou2020data}.) 
In our numerical investigation, we use a Reynolds number $Re=1000$ and four ROM basis functions (i.e., $r=4$).
Details of the computational setting can be found in~\cite{mou2021data}.

In Fig.~\ref{fig:d2-vms-rom-nse}, we present plots of the $L^2$ norm of the error, $\| \bu^{FOM}-\bu^{ROM} \|_{L^2}$, at $t=10$ for G-ROM~\eqref{eqn:g-rom} (top), I-ROM~\eqref{eqn:i-rom} (middle), and D2-VMS-ROM~\eqref{eqn:d2-vms-rom} (bottom).
We note that the G-ROM error is relatively large, whereas the D2-VMS-ROM error is almost negligible.
These two plots clearly show that adding the data-driven closure model to the classical G-ROM (i.e., using the D2-VMS-ROM) significantly increases the G-ROM accuracy.
Although the I-ROM cannot be used in practical computations (since it is not closed), we included I-ROM results for comparison purposes.
The plots in Fig.~\ref{fig:d2-vms-rom-nse} show that the D2-VMS-ROM is not only more accurate than the standard G-ROM, but it is almost as accurate as the I-ROM (which includes an ideal closure model).
Thus, for this test problem, the D2-VMS-ROM error almost reaches the theoretical lower bound given by the I-ROM error.
Overall, Fig.~\ref{fig:d2-vms-rom-nse} clearly shows that closure models can significantly increase the ROM accuracy in under-resolved simulations.

\section{Mathematical Foundations of ROM Closures}
    \label{sec:mathematics}

In Sections~\ref{sec:closure-problem-closure-model} and \ref{sec:d2-vms-rom} we discussed the {\it mathematical modeling} of ROM closures.
In Section~\ref{sec:numerical-results}, we discussed the {\it numerical simulation} of ROM closures.
The following is a natural question:

\begin{question}{Q7}
	What can we prove about ROM closures?
\end{question}

The answer to Q7 is simple:
\begin{important}{A7}
	Not so much.
	Yet.
\end{important}

In this section, we briefly summarize some relevant theoretical aspects associated with ROM closure modeling. 
Compared with the analysis of classical numerical schemes~\cite{BIL05,john2016finite,rebollo2014mathematical}, the theoretical foundations for ROM closures are much less developed.
We emphasize, however, that recently there have been significant advancements in this exciting and important research area. 

The theoretical investigations of ROM closure modeling generally aim at proving error bounds for ROM closures of the form
\begin{eqnarray}
    \| \bu^{FOM} - \bu^{ROM} \| 
    \leq C \left( 
                \text{space error} 
                + \text{time error}
                + \text{ROM error}
            \right) , 
    \label{eqn:math-1}                        
\end{eqnarray}
where $\bu^{FOM}$ is the FOM solution, $\bu^{ROM}$ is the ROM solution, $\| \cdot \|$ is a given norm, the space error is the error that results from the spatial approximation, the time error is the error that results from the time approximation, the ROM error is the error that results from the ROM approximation, and $C$ is a generic constant that does not depend on the discretization parameters.
We note that the first two terms on the right-hand side of~\eqref{eqn:math-1} appear in error bounds for classical numerical discretizations, e.g., the FEM~\cite{john2016finite}.
The third term, however, does not appear in these bounds.

The main purpose of the error bound~\eqref{eqn:math-1} is to show the {\it convergence} of the ROM solution to the FOM solution.
For example, as the spatial mesh size and the time step go to zero, the space error and time error in~\eqref{eqn:math-1}, respectively, are expected to go to zero (at a rate that depends on the particular spatial and time discretizations used).
Furthermore, as the number of ROM basis functions goes to the rank of the snapshot matrix, the ROM error in~\eqref{eqn:math-1} is also expected to go to zero.
Thus, as the right-hand side of~\eqref{eqn:math-1} goes to zero, so does the error on the left-hand side of~\eqref{eqn:math-1}, which proves the convergence of the ROM solution to the FOM solution.

For the G-ROM~\eqref{eqn:g-rom}, the numerical analysis started two decades ago with the pioneering work of Kunisch and Volkwein, who proved the first error bounds for the POD of parabolic equations, e.g., the heat equation~\cite{KV01} and the Navier-Stokes equations~\cite{KV02}. 
More than a decade later, Singler improved Kunisch and Volkwein's results, by proving sharper error bounds~\cite{singler2014new}. 
Recently, optimal pointwise in time error bounds were proved in~\cite{koc2021optimal}.
These results finally bring the G-ROM numerical analysis to a level comparable to (although not as developed as) the level of the numerical analysis of the FEM. 

For the ROM closure models, the numerical analysis is relatively scarce.
The numerical analysis for ROM closures aims at proving a modified form of the G-ROM error bound~\eqref{eqn:math-1}:
\begin{eqnarray}
    \hspace*{-0.6cm}
    \| \bu^{FOM} - \bu^{ROM} \| 
    \leq C \left( 
                \text{space error} 
                + \text{time error}
                + \text{ROM error}
                + \text{closure error}
            \right), 
    \label{eqn:math-2}                        
\end{eqnarray}
where the closure error is the error that results from the approximation of the closure term $\btau^{I-ROM}$ in the I-ROM~\eqref{eqn:i-rom} with a closure model.

As mentioned in~\cite{ahmed2021closures}, the first numerical analysis of ROM closures was performed in~\cite{borggaard2011artificial}, where error bounds for the time discretization of  the Smagorinsky model (i.e., a ROM closure model developed on phenomenological arguments) were proven.
Error bounds for the time and space discretizations of the
Smagorinsky model were later proven in~\cite{rebollo2017certified} in an RBM context.
Error bounds for VMS closure models were proved in~\cite{eroglu2017modular,iliescu2013variational,iliescu2014variational,roop2013proper} (see also~\cite{azaiez2021cure,rubino2020numerical} for related work).
Finally, error bounds for the D2-VMS-ROM~\eqref{eqn:d2-vms-rom} were proved in~\cite{koc2021verifiability} (see also~\cite{koc2019commutation} for related work).

\bigskip

\section{Conclusions and Outlook}
	\label{sec:conclusions}

In this paper, we presented a brief tutorial for reduced order model (ROM) closures. In the first part of our tutorial, we motivated the ROM closures.
We note that ROM closure modeling is often misunderstood in the ROM community.
Thus, we started our tutorial by explaining the need for ROM closure modeling (i.e., the ROM closure problem) in realistic applications, and then we carefully described the ROM closure model.
Specifically, we first outlined the main steps used to construct the Galerkin ROM (G-ROM), which is based on leveraging a data-driven basis in the classical Galerkin framework.
Next, we noted that, although G-ROM can decrease the computational cost of standard numerical discretizations by orders of magnitude, it yields inaccurate results in under-resolved ROM simulations, i.e., when the number of basis functions is not enough to capture the underlying system's dynamics.
To address the G-ROM's inaccuracy in under-resolved simulations, we introduced the ROM closure model.
We motivated the need for ROM closure by presenting a mathematical extension of the classical Galerkin framework to include not only the space of resolved scales, but also the space of unresolved scales.
In this extended variational multiscale framework, we showed that the correct ROM dynamics include an additional term (i.e., the closure term), which represents the effect of the unresolved scales.
Furthermore, we showed that this mathematical framework, which we named the ideal ROM (I-ROM),  
yields numerical results that are significantly more accurate than the G-ROM results.
Thus, we concluded that a ROM closure model, which is a practical model for the I-ROM closure term, should be added to the G-ROM to increase its accuracy in realistic, under-resolved simulations. 

In the second part of our tutorial, we outlined the main steps in the construction of ROM closure models.
To simplify our presentation, we focused on one particular type of ROM closure modeling, i.e., data-driven modeling.
Furthermore, we illustrated this construction for one specific data-driven ROM closure model, i.e., the data-driven variational multiscale ROM (D2-VMS-ROM).
In our construction, we started with the closure term in the I-ROM and we simply posed the closure problem as leveraging the available FOM data to find the ``best'' ROM closure model. 
To this end, we first postulated a model form for the ROM closure model. 
Then, we solved a least squares problem to find the parameters in the model form that yield the ROM closure model that is the closest to the ideal ROM closure model.
Finally, we also included numerical results for the two-dimensional flow past a circular cylinder, which showed that the D2-VMS-ROM was significantly more accurate than the standard G-ROM, and almost as accurate as the I-ROM.
These numerical results illustrated the significant benefit of ROM closure modeling in under-resolved simulations.

We hope that this brief tutorial offers a glimpse into the exciting research field of ROM closure modeling, which has witnessed a significant development over the past two decades.
This research area is currently experiencing a dynamic development in several directions.
One of the most active research directions is the use of machine learning tools to construct more accurate and more efficient ROM closure models.
Recently, deep learning models have been shown to be quite effective and computationally efficient in capturing the relationship between resolved and unresolved scales \cite{ahmed2020long}. 
However, these models often need large amounts of training data and their generalization, expressivity, and analysis still remain mostly challenging.

Another important research direction is the development of ROM closures for problems in solid mechanics.
Although most ROM closure modeling has been performed in computational fluid dynamics~\cite{ahmed2021closures}, there has been recent work done in solid mechanics. 
For example, approximations of the mechanical behavior of 
soft tissue showed substantial improvement in accuracy over G-ROM with the addition of ROM closure terms at a modest computational overhead~\cite{snyder2022data}. 
The ability of ROM closure to capture the nonlinearities of soft tissue behavior is especially promising for its application in biomechanics.

Depending on the applications, one 
can also couple ROMs with additional parameterization schemes or surrogate models for some of the unresolved scales in order to recover more dynamical features of the original system, especially when the ROMs are constructed for under-resolved dynamical regimes. 
For instance, in the context of data assimilation, when observations are only available for the (large-scale) low-frequency modes, one can 
design computationally efficient strategies within the conditional Gaussian framework \cite{chen2020learning,chen2021conditional,chen2018conditional} to approximate the dynamics of the high-frequency (unresolved) modes with quantified uncertainties by a suitable dynamical model for the unresolved modes.

Finally, providing mathematical support for ROM closures is also an important research direction.
We note that significant mathematical support has been provided for closures in classical computational fluid dynamics~\cite{BIL05,john2016finite,rebollo2014mathematical}.
For ROM closures, however, only the first steps have been taken
and much more remains to be done.

\bibliographystyle{plain}
\bibliography{traian,honghu}

\end{document}